\documentclass[12pt]{amsart}
\usepackage{latexsym,amssymb,amsrefs}
\usepackage[all]{xy}

\newtheorem{lemma}{Lemma}

\newtheorem{theorem}[lemma]{Theorem}

\theoremstyle{definition}
\newtheorem{example}{Example}

\newcommand{\textmc}[1]{\textsf{#1}}

\renewcommand{\MR}{\ }

\begin{document}

\title[Multicomplexes and spectral sequences]{Multicomplexes and spectral sequences}

\author{David E. Hurtubise}
\address{Department of Mathematics and Statistics\\
         Penn State Altoona\\
         Altoona, PA 16601-3760}
\email{Hurtubise@psu.edu}
\subjclass[2000]{Primary: 18G40; 55T99}

\begin{abstract}
In this note we present some algebraic examples of multicomplexes whose
differentials differ from those in the spectral sequences associated
to the multicomplexes. The motivation for constructing examples  
showing the algebraic distinction between a multicomplex and its 
associated spectral sequence comes from the author's work on Morse-Bott 
homology with A. Banyaga \cite{BanMor}.
\end{abstract}

\maketitle


\section{Introduction}

Let $R$ be a principal ideal domain. A first quadrant \textbf{multicomplex} $X$ is 
a bigraded $R$-module $\{X_{p,q}\}_{p,q \in \mathbb{Z}_+}$ with
differentials
$$
\textmc{d}_i:X_{p,q} \rightarrow X_{p-i,q+i-1} \quad \text{ for all } i = 0,1,\ldots
$$
that satisfy
$$
\sum_{i+j = n} \textmc{d}_i \textmc{d}_j = 0 \quad \text{ for all }n.
$$
A first quadrant multicomplex such that $\textmc{d}_i = 0$ for all $i \geq 2$ is called
a \textbf{double complex} (or a \textbf{bicomplex}). For the basic properties of
multicomplexes we refer the reader to \cite{BoaCon} and \cite{MeyAcy}.

\smallskip
An $E^k$ first quadrant \textbf{spectral sequence} is a sequence of bigraded $R$-modules
$\{E^r_{s,t}\}_{s,t \in \mathbb{Z}_+}$ with differentials
$$
d^r:E^r_{s,t} \rightarrow E^r_{s-r,t+r-1}
$$
such that for all $r \geq k$ there is a given isomorphism $H(E^r) \approx E^{r+1}$
(see Section \ref{ss}). Every first quadrant multicomplex determines 
an $E^0$ first quadrant spectral sequence. However, not every first 
quadrant spectral sequence comes from a multicomplex.

\smallskip
Moreover, the differentials $d^r$ in the spectral sequence associated to a first quadrant 
multicomplex are in general different from the homomorphisms induced by the
differentials $\textmc{d}_i$ in the multicomplex. (Note that using the term ``differential''
to describe the homomorphisms $\textmc{d}_i$ in a multicomplex is misleading since there is no 
guarantee that $(\textmc{d}_i)^2$ is zero.) The purpose of this note is to demonstrate this
distinction by presenting explicit algebraic examples of multicomplexes where the differential
$d^r$ in the associated spectral sequence is different than the homomorphism induced by $\textmc{d}_r$.

\medskip
The motivation for constructing examples that show the distinction between 
a multicomplex and its associated spectral sequence comes from the author's work
on Morse-Bott homology with A. Banyaga and the discovery that the
Morse-Bott-Smale chain complex is in fact a multicomplex.  For more details see
the introduction to \cite{BanMor}.


\section{The spectral sequence associated to a filtered chain complex}\label{ss}

In this section we clarify the meaning (and the bigrading) of the isomorphism
$H(E^r) \approx E^{r+1}$, and we recall the definition of the differentials
$d^r:E^r_{s,t} \rightarrow E^r_{s-r,t+r-1}$ in an $E^k$ spectral sequence coming from
a filtered chain complex. This section follows Chapter 9 of \cite{SpaAlg}. 

\smallskip
An $E^k$ spectral sequence consists of a sequence of bigraded modules $\{E^r_{s,t}\}$ over a
principal ideal domain $R$ for $r \geq k$, with differentials
$d^r:E^r_{s,t} \rightarrow E^r_{s-r,t+r-1}$ that satisfy $(d^r)^2 = 0$. If we define
\begin{eqnarray*}
\bar{Z}_{s,t}^r & = & \text{ker}(d^r:E_{s,t}^r \rightarrow E_{s-r,t+r-1}^r)\\
\bar{B}_{s,t}^r & = & \text{im}(d^r:E_{s+r,t-r+1}^r \rightarrow E_{s,t}^r)
\end{eqnarray*}
then $\bar{B}^r_{s,t} \subseteq \bar{Z}^r_{s,t}$, and by definition there is a given
isomorphism
$
E^{r+1}_{s,t} \approx \bar{Z}^r_{s,t}/\bar{B}^r_{s,t}.
$

\smallskip
Let $(C_\ast,\partial)$ be a filtered chain complex that is bounded below
by $s = 0$. That is, suppose that we have a filtration
$$
F_0 C_\ast \subset \cdots \subset F_{s-1}C_\ast \subset F_s C_\ast \subset F_{s+1} 
C_\ast \subset \cdots
$$
where $F_sC_\ast$ is a chain subcomplex of $C_\ast$ for all $s$, i.e. 
$\partial(F_sC_{s+t}) \subseteq F_sC_{s+t-1}$ for all $t$.  
The grading $s$ is called the \textbf{filtered degree}, the grading $t$ is called the
\textbf{complementary degree}, and the sum $s+t$ is called the \textbf{total degree}.
The filtration is said to be \textbf{convergent} if $\cap_s F_sC_\ast = 0$ and $\cup_s F_sC_\ast = C_\ast$.
Define
\begin{eqnarray*}
Z_{s,t}^r & = & \{ c \in F_sC_{s+t}|\ \partial c \in F_{s-r} C_{s+t-1}\}\\
Z_{s,t}^\infty & = & \{ c \in F_sC_{s+t}|\ \partial c = 0\}.
\end{eqnarray*}
The bigraded $R$-modules in the spectral sequence associated to the filtration
are defined to be
\begin{eqnarray*}
E_{s,t}^r & = & Z_{s,t}^r\left/\left(Z^{r-1}_{s-1,t+1} + \partial Z^{r-1}_{s+r-1,t-r+2} 
 \right)\right.\\
E_{s,t}^\infty & = & Z_{s,t}^\infty\left/\left(Z^\infty_{s-1,t+1} + (\partial C_{s+t+1} \cap
 F_sC_{s+t}) \right)\right.,
\end{eqnarray*}
where $A+B$ denotes the free abelian group generated by the elements of $A$ and $B$,
and the differential $d^r:E^r_{s,t} \rightarrow E^r_{s-r,t+r-1}$ is defined
by the following diagram.
$$
\xymatrix{
Z^r_{s,t} \ar[r]^(.45){\partial} \ar[d] & Z^r_{s-r,t+r-1} \ar[d]\\
Z_{s,t}^r\left/\left(Z^{r-1}_{s-1,t+1} + \partial Z^{r-1}_{s+r-1,t-r+2}\right)\right. 
 \ar[r]^(.48){d^r} &
Z_{s-r,t+r-1}^r\left/\left(Z^{r-1}_{s-r-1,t+r} + \partial Z^{r-1}_{s-1,t+1}\right)\right.
}
$$
The $R$-module $E^{r}_{s,t}$ is isomorphic to $\bar{Z}^{r-1}_{s,t}/\bar{B}^{r-1}_{s,t}$ 
via an isomorphism given by the Noether Isomorphism Theorem.

\smallskip
For a proof of the following theorem see Section 9.1 of \cite{SpaAlg}.

\begin{theorem}\label{specfilt}
If the filtration on the chain complex $(C_\ast,\partial)$ is convergent
and bounded below, then the above spectral sequence converges to the 
bigraded $R$-module $GH_\ast(C_\ast,\partial)$ associated to the filtration 
$F_sH_\ast(C_\ast,\partial) \equiv \text{im}[H_\ast(F_sC_\ast,\partial)
\rightarrow H_\ast(C_\ast,\partial)]$. That is,
$$
E^\infty_{s,t} \approx \bigcap_r Z^r_{s,t} \left/ \bigcup_r \left(Z^{r-1}_{s-1,t+1} + 
\partial Z^{r-1}_{s+r-1,t-r+2} \right) \right. \approx GH_\ast(C_\ast,\partial)_{s,t}
$$
where $GH_\ast(C_\ast,\partial)_{s,t} \equiv 
F_sH_{s+t}(C_\ast,\partial)/F_{s-1}H_{s+t}(C_\ast,\partial)$.
\end{theorem}

%


\section{The spectral sequence associated to a multicomplex}

A first quadrant multicomplex $(\{X_{p,q}\}_{p,q \in \mathbb{Z}_+},
\{\textmc{d}_i\}_{i \in \mathbb{Z}_+})$ can be \textbf{assembled} to form a filtered 
chain complex $((CX)_\ast,\partial)$ by summing along the diagonals. That is, 
suppose that we are given a bigraded $R$-module $\{X_{p,q}\}_{p,q \in \mathbb{Z}_+}$ 
and homomorphisms 
$$
\textmc{d}_i:X_{p,q} \rightarrow X_{p-i,q+i-1} \quad \text{ for all } i = 0,1,\ldots
$$
that satisfy
$$
\sum_{i+j = n} \textmc{d}_i \textmc{d}_j = 0 \quad \text{ for all }n.
$$

$$
\xymatrix{
\vdots & \vdots & \vdots & \vdots \\
X_{0,3} \ar[d]^(.4){\textmc{d}_0} & X_{1,3} \ar[d]^(.4){\textmc{d}_0} \ar[l]_{\textmc{d}_1} & X_{2,3} \ar[d]^(.4){\textmc{d}_0} \ar[l]_{\textmc{d}_1} & X_{3,3} \ar[d]^(.4){\textmc{d}_0}
\ar[l]_{\textmc{d}_1} & \cdots\\
X_{0,2} \ar[d]^(.4){\textmc{d}_0} & X_{1,2} \ar[d]^(.4){\textmc{d}_0} \ar[l]_{\textmc{d}_1} & X_{2,2} \ar[d]^(.4){\textmc{d}_0} \ar[l]_{\textmc{d}_1} \ar[llu]|(.37){\textmc{d}_2} |!{[ul];[l]}\hole & 
X_{3,2} \ar[d]^(.4){\textmc{d}_0} \ar[l]_{\textmc{d}_1} 
\ar[llu]|(.37){\textmc{d}_2} |!{[ul];[l]}\hole & \cdots\\
X_{0,1} \ar[d]^(.4){\textmc{d}_0} & X_{1,1} \ar[d]^(.4){\textmc{d}_0} \ar[l]_{\textmc{d}_1} & X_{2,1} \ar[d]^(.4){\textmc{d}_0} \ar[l]_{\textmc{d}_1} \ar[llu]|(.37){\textmc{d}_2} |!{[ul];[l]}\hole & 
X_{3,1} \ar[d]^(.4){\textmc{d}_0} \ar[l]_{\textmc{d}_1} \ar[llu]|(.37){\textmc{d}_2} |!{[ul];[l]}\hole 
  \ar@{.>}[llluu]_>(.2){\textmc{d}_3} & \cdots \\
X_{0,0}  & X_{1,0} \ar[l]_{\textmc{d}_1} & X_{2,0} \ar[l]_{\textmc{d}_1} \ar[llu]|(.37){\textmc{d}_2} |!{[ul];[l]}\hole & X_{3,0} \ar[l]_{\textmc{d}_1} \ar[llu]|(.37){\textmc{d}_2} |!{[ul];[l]}\hole \ar@{.>}[llluu]_>(.2){\textmc{d}_3} & \cdots
}
$$

If we define
$$
(CX)_{n} \equiv \bigoplus_{p+q = n} X_{p,q}
$$
and $\partial_n = \textmc{d}_0 \oplus \cdots \oplus \textmc{d}_n$ for all $n \in \mathbb{Z}_+$,
then the above relations imply that $\partial_{n} \circ \partial_{n+1} = 0$. 
$$
\xymatrix{
\cdots & X_{3,0} \ar[r]^(.53){\textmc{d}_0}  \ar[dr]^{\textmc{d}_1} 
\ar[ddr]|(.33){\textmc{d}_2} |!{[d];[dr]}\hole \ar@{.>}[dddr]|(.25){\textmc{d}_3} & 0 \\
\cdots & X_{2,1} \ar@{}[u]|{\oplus} \ar[r]^(.53){\textmc{d}_0} \ar[dr]^{\textmc{d}_1} 
\ar[ddr]|(.33){\textmc{d}_2} |!{[d];[dr]}\hole & X_{2,0} \ar@{}[u]|{\oplus} \ar[r]^(.53){\textmc{d}_0} 
\ar[dr]^{\textmc{d}_1} \ar[ddr]|(.33){\textmc{d}_2} |!{[d];[dr]}\hole & 0 & & \\
\cdots & X_{1,2} \ar@{}[u]|{\oplus} \ar[r]^(.53){\textmc{d}_0} \ar[dr]^{\textmc{d}_1} & X_{1,1} 
\ar@{}[u]|{\oplus} \ar[r]^(.53){\textmc{d}_0} \ar[dr]^{\textmc{d}_1} & X_{1,0} \ar@{}[u]|{\oplus} 
\ar[r]^{\textmc{d}_0} \ar[dr]^{\textmc{d}_1} & 0 & \\
\cdots & X_{0,3} \ar@{}[u]|{\oplus} \ar[r]^(.53){\textmc{d}_0} & X_{0,2} \ar@{}[u]|{\oplus} 
\ar[r]^(.53){\textmc{d}_0} & X_{0,1} \ar@{}[u]|{\oplus} \ar[r]^(.53){\textmc{d}_0} & X_{0,0} \ar@{}[u]|{\oplus}
\ar[r]^-{\textmc{d}_0} & 0\\
\cdots &  {(CX)}_3 \ar@{}[u]|{\|} \ar[r]^{\partial_3} & {(CX)}_2 \ar@{}[u]|{\|} 
\ar[r]^{\partial_2} & {(CX)}_1 \ar@{}[u]|{\|} \ar[r]^{\partial_1} & {(CX)}_0 \ar@{}[u]|{\|}
\ar[r]^-{\partial_0} & 0 \ar@{}[u]|{\|}
}
$$
Moreover, the chain complex $((CX)_\ast,\partial)$ has an obvious  filtration given by
$$
F_s(CX)_n\ \equiv \bigoplus_{\stackrel{\scriptstyle p+q = n}{p\leq s}} X_{p,q}.
$$
Note that the restriction $q \leq s$ determines a second filtration on a double complex,
but it does not determine a filtration on a general multicomplex.

\medskip
The bigraded module associated to the above filtration is
$$
G((CX)_\ast)_{s,t} = F_s(CX)_{s+t}/F_{s-1}(CX)_{s+t} \approx X_{s,t}
$$
for all $s,t \in \mathbb{Z}_+$, and the $E^1$ term of the associated spectral sequence
is given by
$$
E_{s,t}^1 = Z_{s,t}^1\left/\left(Z^{0}_{s-1,t+1} + \partial Z^{0}_{s,t+1}\right)\right.
$$
where
$$
Z^1_{s,t} = \{c \in F_s(CX)_{s+t} |\ \partial c \in F_{s-1}(CX)_{s+t-1}\} 
$$
and
$$
Z^0_{s,t} = \{c \in F_s(CX)_{s+t} |\ \partial c \in F_{s}(CX)_{s+t-1}\} = F_s(CX)_{s+t}.
$$
The group $Z^1_{s,t}/Z^0_{s-1,t+1}$ is the group of $(s+t)$-cycles in the quotient chain complex
$F_s(CX)_\ast/F_{s-1}(CX)_\ast$, and the group of $(s+t)$-boundaries in
$F_s(CX)_\ast/F_{s-1}(CX)_\ast$ is
$$
\partial Z^0_{s,t+1} \left/ (Z^0_{s-1,t+1} \cap \partial Z^0_{s,t+1}) \right. \approx
\left(Z^0_{s-1,t+1} + \partial Z^0_{s,t+1} \right) \left/ Z^0_{s-1,t+1}\right.
$$
where the isomorphism is given by the Noether Isomorphism Theorem. Therefore,
\begin{eqnarray*}
E^1_{s,t} & = & Z_{s,t}^1\left/\left(Z^{0}_{s-1,t+1} + \partial Z^{0}_{s,t+1}\right)\right.\\
& \approx & \frac{Z_{s,t}^1/Z^0_{s-1,t+1}}
  {\left(Z^{0}_{s-1,t+1} + \partial Z^{0}_{s,t+1}\right)/Z^0_{s-1,t+1}}\\
& \approx & \frac{Z_{s,t}^1/Z^0_{s-1,t+1}}
  {\partial Z^0_{s,t+1} \left/ (Z^0_{s-1,t+1} \cap \partial Z^0_{s,t+1}) \right.}
\end{eqnarray*}
and we see that $E^1_{s,t} \approx H_{s+t}(X_{s,\ast},\textmc{d}_0)$ where $(X_{s,\ast},\textmc{d}_0)$
denotes the chain complex
$$
\xymatrix{
\cdots  \ar[r]^{\textmc{d}_0} &  X_{s,3} \ar[r]^{\textmc{d}_0} & X_{s,2} \ar[r]^{\textmc{d}_0} & X_{s,1} \ar[r]^{\textmc{d}_0} & 
X_{s,0} \ar[r]^{\textmc{d}_0} & 0.
}
$$

The differential $d^1$ on the $E^1$ term of the spectral sequence is defined by the
diagram
$$
\xymatrix{
Z^1_{s,t} \ar[r]^(.47){\partial} \ar[d] & Z^1_{s-1,t} \ar[d]\\
Z_{s,t}^1\left/\left(Z^{0}_{s-1,t+1} + \partial Z^{0}_{s,t+1}\right)\right. 
 \ar[r]^(.47){d^1} &
Z^1_{s-1,t}\left/\left(Z^{0}_{s-2,t+1} + \partial Z^{0}_{s-1,t+1}\right)\right.
}
$$
and it is natural to ask whether or not there is a connection between the differential
$d^1:E^1_{s,t} \rightarrow E^1_{s-1,t}$ in the spectral sequence and the homomorphism 
$\textmc{d}_1:X_{s,t} \rightarrow X_{s-1,t}$ in the multicomplex.  

It is an easy exercise to show that the relations 
\begin{eqnarray*}
\textmc{d}_0 \textmc{d}_1 + \textmc{d}_1 \textmc{d}_0 & = & 0\\ 
\textmc{d}_0\textmc{d}_2 + \textmc{d}_1 \textmc{d}_1 + \textmc{d}_2 \textmc{d}_0 & = & 0
\end{eqnarray*}
imply that the homomorphism $\textmc{d}_1$ induces a differential $\bar{\textmc{d}}_1:E^1_{s,t}
\rightarrow E^1_{s-1,t}$, i.e. $(\bar{\textmc{d}}_1)^2 = 0$.  Moreover, one can show that the
differential $\bar{\textmc{d}}_1$ coincides with the differential $d^1:E^1_{s,t} \rightarrow 
E^1_{s-1,t}$. This is a standard fact for double complexes, and the proof for double 
complexes carries over to multicomplexes (see for instance Section 14 of \cite{BotDif},
or Section 3.2.1 of \cite{McCAus}).  Thus, we have the following.

\begin{theorem}\label{r=1}
Let $(\{X_{p,q}\}_{p,q \in \mathbb{Z}_+},\{\textmc{d}_i\}_{i \in \mathbb{Z}_+})$ be a first 
quadrant multicomplex and $((CX)_\ast,\partial)$ the associated assembled chain complex.
Then the $E^1$ term of the spectral sequence associated to the filtration of $(CX)_\ast$
determined by the restriction $p \leq s$ is given by $E^1_{s,t} \approx 
H_{s+t}(X_{s,\ast},\textmc{d}_0)$ where $(X_{s,\ast},\textmc{d}_0)$ denotes the following chain complex.
$$
\xymatrix{
\cdots  \ar[r]^{\textmc{d}_0} &  X_{s,3} \ar[r]^{\textmc{d}_0} & X_{s,2} \ar[r]^{\textmc{d}_0} & X_{s,1} \ar[r]^{\textmc{d}_0} & 
X_{s,0} \ar[r]^{\textmc{d}_0} & 0
}
$$
Moreover, the $d^1$ differential on the $E^1$ term of the spectral sequence is induced
from the homomorphism $\textmc{d}_1$ in the multicomplex.
\end{theorem}


\section{Multicomplexes where $d^r \neq \textbf{d}_r$}

Theorem \ref{r=1} should sound familiar to anyone acquainted with  
double complexes. However, the examples in this section
show that  Theorem \ref{r=1} does \textbf{not} generalize to 
the higher differentials in the spectral sequence associated
to a multicomplex. In fact, the pattern suggested by Theorem \ref{r=1}
breaks down when $r=2$. That is, the differential $d^r$ in the spectral
sequence is \textbf{not} necessarily induced from the homomorphism 
$\textmc{d}_r$ when $r \geq 2$. To paraphrase Section 11 of \cite{BoaCon}, 
when $r \geq 2$ the differential $d^r$ is induced from $\textmc{d}_r$ only 
on those classes which contain elements $x$ such that $\textmc{d}_i(x) = 0$ 
for all $i < r$ ``which rarely happens''.


\begin{example}[A double complex with $d^2 \neq 0$]\label{doublecomplex}
\end{example}

It is well known that the spectral sequence associated to a double complex
does not necessarily degenerate at $E^2$.  That is, there is no guarantee
that $d^r=0$ for $r \geq 2$. This first example is a small algebraic 
example that demonstrates this phenomena.  

\medskip
Consider the following first quadrant double complex
$$
\xymatrix{
0 \ar[d]^(.4){\textmc{d}_0} & 0 \ar[d]^(.4){\textmc{d}_0} \ar[l]_(.45){\textmc{d}_1} & 0 \ar[d]^(.4){\textmc{d}_0} \ar[l]_(.45){\textmc{d}_1} \\
<x_{0,1}> \ar[d]^(.4){\textmc{d}_0} & <x_{1,1}> \ar[d]^(.4){\textmc{d}_0} \ar[l]_(.45){\textmc{d}_1} & 0 \ar[d]^(.4){\textmc{d}_0} \ar[l]_(.45){\textmc{d}_1} \\
0  & <x_{1,0}> \ar[l]_(.45){\textmc{d}_1} & <x_{2,0}> \ar[l]_(.45){\textmc{d}_1}
}
$$
where $<x_{p,q}>$ denotes the free abelian group generated by $x_{p,q}$, the
groups $X_{p,q} = 0$ for $p+q > 2$, and the homomorphisms $\textmc{d}_0$ and $\textmc{d}_1$ satisfy the
following: $\textmc{d}_0(x_{1,1}) = x_{1,0}$, $\textmc{d}_1(x_{1,1}) = x_{0,1}$, and $\textmc{d}_1(x_{2,0}) = x_{1,0}$.
The conditions $(\textmc{d}_0)^2 = (\textmc{d}_1)^2 = 0$ and $\textmc{d}_0 \textmc{d}_1 + \textmc{d}_1 \textmc{d}_0 = 0$
are satisfied trivially, and the assembled chain complex associated
to this double complex is as follows.
$$
\xymatrix{
\cdots & 0 \ar[r]^(.4){\textmc{d}_0}  \ar[dr]^{\textmc{d}_1} & 0 \\
\cdots & 0 \ar@{}[u]|{\oplus} \ar[r]^(.4){\textmc{d}_0} \ar[dr]^{\textmc{d}_1} & <x_{2,0}> \ar@{}[u]|{\oplus} \ar[r]^{\textmc{d}_0} 
\ar[dr]^{\textmc{d}_1} & 0 & & \\
\cdots & 0 \ar@{}[u]|{\oplus} \ar[r]^(.4){\textmc{d}_0} \ar[dr]^{\textmc{d}_1} & <x_{1,1}> 
\ar@{}[u]|{\oplus} \ar[r]^{\textmc{d}_0} \ar[dr]^{\textmc{d}_1} & <x_{1,0}> \ar@{}[u]|{\oplus} 
\ar[r]^(.63){\textmc{d}_0} \ar[dr]^{\textmc{d}_1} & 0 & \\
\cdots & 0 \ar@{}[u]|{\oplus} \ar[r]^(.4){\textmc{d}_0} & 0 \ar@{}[u]|{\oplus} 
\ar[r]^{\textmc{d}_0} & <x_{0,1}> \ar@{}[u]|{\oplus} \ar[r]^(.63){\textmc{d}_0} & 0 \ar@{}[u]|{\oplus}
\ar[r]^-{\textmc{d}_0} & 0\\
\cdots & 0 \ar@{}[u]|{\|} \ar[r]^(.4){\partial_3} & {(CX)}_2 \ar@{}[u]|{\|} 
\ar[r]^{\partial_2} & {(CX)}_1 \ar@{}[u]|{\|} \ar[r]^(.63){\partial_1} & 0 \ar@{}[u]|{\|}
\ar[r]^-{\partial_0} & 0 \ar@{}[u]|{\|}
}
$$

The homology $H_n((CX)_\ast,\partial)$ of the assembled chain 
complex is trivial for all $n \in \mathbb{Z}_+$ because the kernel of $\partial_2$
is trivial and both $x_{0,1}$ and $x_{1,0}$ are in the image of $\partial_2 = \textmc{d}_0 + \textmc{d}_1$:
\begin{eqnarray*}
\partial_2 (x_{2,0}) & = & x_{1,0}\\
\partial_2 (x_{1,1} - x_{2,0}) & = & x_{0,1}.
\end{eqnarray*}
However, the $E^1$ term of the associated spectral sequence is
$$
\xymatrix{
0         & 0 \ar[l]_(.3){\textmc{d}_1} & 0 \ar[l]_(.58){\textmc{d}_1}  \\
<x_{0,1}> & 0 \ar[l]_(.3){\textmc{d}_1} & 0 \ar[l]_(.58){\textmc{d}_1}  \\
0         & 0 \ar[l]_(.3){\textmc{d}_1} & <x_{2,0}> \ar[l]_(.58){\textmc{d}_1} 
}
$$
where $E^1_{s,t} = 0$ for all $s+t > 2$, and the $E^2$ term is isomorphic to
the $E^1$ term. Since $H_n((CX)_\ast,\partial) = 0$ for all $n \in \mathbb{Z}_+$,
Theorem \ref{specfilt} implies that the differential $d^2$ in the spectral sequence
must be nonzero, even though the homomorphism $\textmc{d}_2$ in the multicomplex is zero,
i.e. $d^2\neq 0$ is \textbf{not} induced from $\textmc{d}_2 = 0$.

To compute the differential $d^2:E^2_{2,0} \rightarrow E^2_{0,1}$ we consider 
the following diagram
$$
\xymatrix{
Z^2_{2,0} \ar[r]^(.5){\partial} \ar[d] & Z^2_{0,1} \ar[d]\\
Z^2_{2,0}\left/\left(Z^1_{1,1} + \partial Z^1_{3,0}\right)\right. 
 \ar[r]^(.5){d^2} &
Z^2_{0,1}\left/\left(Z^1_{-1,2} + \partial Z^1_{1,1}\right)\right.
}
$$
where 
$$
Z^2_{2,0} = <x_{1,1} - x_{2,0}>, \quad Z^2_{0,1} = < x_{0,1}>
$$
and $Z^1_{1,1} = Z^1_{3,0} = Z^1_{-1,2} = 0$. Since $\partial_2(x_{1,1} - x_{2,0}) = x_{0,1}$,
we see that the $E^3$ term of the spectral sequence is trivial, and we have
verified Theorem \ref{specfilt} in this example.  

There is one more subtle point to note in this example: although $E^1_{2,0}$ and $E^2_{2,0}$
are isomorphic, they have different generators. That is,
$$
E^1_{2,0} = Z^1_{2,0}\left/\left(Z^0_{1,1} + \partial Z^0_{2,1}\right)\right. \approx\ <x_{2,0}>
$$
whereas
$$
E^2_{2,0} = Z^2_{2,0}\left/\left(Z^1_{1,1} + \partial Z^1_{3,0}\right)\right. \approx\ <x_{1,1} - x_{2,0}>.
$$
This is consistent with the definition of a spectral sequence which states that 
``there is a given isomorphism $H(E^r) \approx E^{r+1}$''
\cite{SpaAlg}.


\begin{example}[A double complex with some $d^r \neq 0$ for $r$ arbitrarily large]\label{double_r}
\end{example}

The preceeding example can be generalized to produce a double complex such
that a differential $d^r$ in the associated spectral sequence is
nonzero for $r$ arbitrarily large. To see this pick any $r\in \mathbb{Z}_+$
with $r \geq 2$, and consider the following first quadrant double complex
$$
\hspace{-.27 in}
\xymatrix{
0 \ar[d]^{\textmc{d}_0} & 0 \ar[d]^{\textmc{d}_0} \ar[l]_{\textmc{d}_1} & 0 \ar[d]^{\textmc{d}_0} \ar[l]_{\textmc{d}_1} & 0 \ar[d]^{\textmc{d}_0} \ar[l]_{\textmc{d}_1} & \cdots \ar[l]_{\textmc{d}_1}
& 0 \ar[d]^{\textmc{d}_0} \ar[l]_(.43){\textmc{d}_1} \\
<x_{0,r-1}> \ar[d]^{\textmc{d}_0} & <x_{1,r-1}> \ar[d]^{\textmc{d}_0} \ar[l]_{\textmc{d}_1} & 0 \ar[d]^{\textmc{d}_0} \ar[l]_{\textmc{d}_1} & 0 \ar[l]_{\textmc{d}_1} \ar[d]^{\textmc{d}_0}& 
  \cdots \ar[l]_{\textmc{d}_1} & 0 \ar[l]_(.43){\textmc{d}_1} \ar[d]^{\textmc{d}_0} \\
0 \ar[d]^{\textmc{d}_0} & <x_{1,r-2}> \ar[l]_{\textmc{d}_1} \ar[d]^{\textmc{d}_0} & <x_{2,r-2}> \ar[d]^{\textmc{d}_0} \ar[l]_{\textmc{d}_1} & 0 \ar[d]^{\textmc{d}_0} 
  \ar[l]_{\textmc{d}_1} & \cdots  \ar[l]_{\textmc{d}_1} & 0 \ar[l]_(.43){\textmc{d}_1} \ar[d]^{\textmc{d}_0} \\
\vdots \ar[d]^{\textmc{d}_0} & \vdots & \vdots  \ar[d]^{\textmc{d}_0} & \vdots  \ar[d]^{\textmc{d}_0} & \vdots  \ar[d]^{\textmc{d}_0} & \vdots  \ar[d]^{\textmc{d}_0} \\
0 \ar[d]^{\textmc{d}_0} & \cdots  \ar[l]_{\textmc{d}_1} & 0 \ar[d]^{\textmc{d}_0} \ar[l]_{\textmc{d}_1} & <x_{r-2,1}> \ar[d]^{\textmc{d}_0} \ar[l]_{\textmc{d}_1} & <x_{r-1,1}> \ar[d]^{\textmc{d}_0}
   \ar[l]_{\textmc{d}_1} & 0 \ar[d]^{\textmc{d}_0} \ar[l]_(.43){\textmc{d}_1}\\
0  & \cdots  \ar[l]_{\textmc{d}_1} & 0 \ar[l]_{\textmc{d}_1} & 0 \ar[l]_{\textmc{d}_1} & <x_{r-1,0}> \ar[l]_{\textmc{d}_1} & <x_{r,0}> \ar[l]_(.43){\textmc{d}_1}
}
$$
where the groups $X_{p,q} = 0$ for $p+q > r$ and the homomorphisms $\textmc{d}_0$ and 
$\textmc{d}_1$ satisfy the following for $p+q = r$: $\textmc{d}_0(x_{p,q}) = x_{p,q-1}$ and
$\textmc{d}_1(x_{p,q}) = x_{p-1,q}$. The conditions $(\textmc{d}_0)^2 = (\textmc{d}_1)^2 = 0$
and $\textmc{d}_0 \textmc{d}_1 + \textmc{d}_1 \textmc{d}_0 = 0$ are satisfied trivially,
and the assembled chain complex associated to this double complex is as follows.

$$
\xymatrix{
\cdots & 0 \ar[r]^(.35){\textmc{d}_0}  \ar[dr]^{\textmc{d}_1} & 0 \\
\cdots & 0 \ar@{}[u]|{\oplus} \ar[r]^(.35){\textmc{d}_0} \ar[dr]^{\textmc{d}_1} & <x_{r,0}> \ar@{}[u]|{\oplus} \ar[r]^{\textmc{d}_0} 
\ar[dr]^{\textmc{d}_1} & 0 & & \\
\cdots & 0 \ar@{}[u]|{\oplus} \ar[r]^(.35){\textmc{d}_0} \ar[dr]^{\textmc{d}_1} & <x_{r-1,1}> 
\ar@{}[u]|{\oplus} \ar[r]^{\textmc{d}_0} \ar[dr]^{\textmc{d}_1} & <x_{r-1,0}> \ar@{}[u]|{\oplus} 
\ar[r]^(.65){\textmc{d}_0} \ar[dr]^{\textmc{d}_1} & 0 & \\
 & \vdots & \vdots & \vdots & \vdots \\
\cdots & 0 \ar@{}[u]|{\oplus} \ar[r]^(.35){\textmc{d}_0} \ar[dr]^{\textmc{d}_1} & <x_{1,r-1}> 
\ar@{}[u]|{\oplus} \ar[r]^{\textmc{d}_0} \ar[dr]^{\textmc{d}_1} & <x_{1,r-2}> \ar@{}[u]|{\oplus} 
\ar[r]^(.65){\textmc{d}_0} \ar[dr]^{\textmc{d}_1} & 0 \ar@{}[u]|{\oplus} & \cdots \\
\cdots & 0 \ar@{}[u]|{\oplus} \ar[r]^(.35){\textmc{d}_0} & 0 \ar@{}[u]|{\oplus} 
\ar[r]^{\textmc{d}_0} & <x_{0,r-1}> \ar@{}[u]|{\oplus} \ar[r]^(.65){\textmc{d}_0} & 0 \ar@{}[u]|{\oplus} & \cdots \\
\cdots & 0 \ar@{}[u]|{\|} \ar[r]^(.35){\partial_{r+1}} & {(CX)}_r \ar@{}[u]|{\|} 
\ar[r]^{\partial_r} & {(CX)}_{r-1} \ar@{}[u]|{\|} \ar[r]^(.65){\partial_{r-1}} & 0 \ar@{}[u]|{\|} & \cdots
}
$$

As in the previous example, the homology $H_n((CX)_\ast,\partial)$ of
the assembled chain complex is trivial for all $n \in \mathbb{Z}_+$ 
because the kernel of $\partial_r$ is trivial and all the generators
$x_{0,r-1}, x_{1,r-2}, \ldots ,x_{r-1,0}$ of $(CX)_{r-1}$ are in the 
image of $\partial_r$. The $E^1$ term of the associated spectral sequence
has $E^1_{0,r-1} =\ <x_{0,r-1}>$, $E^1_{r,0} =\ <x_{r,0}>$, and $E^1_{s,t} = 0$
for all other values of $s$ and $t$.
Moreover, $E^1 \approx E^2 \approx \cdots \approx E^{r-1}$.
Once again, Theorem \ref{specfilt} implies that the differential $d^r$ must be nonzero
(even though $\textmc{d}_r = 0$), and the diagram
$$
\xymatrix{
Z^r_{r,0} \ar[r]^(.52){\partial} \ar[d] & Z^r_{0,r-1} \ar[d]\\
Z_{r,0}^r\left/\left(Z^{r-1}_{r-1,1} + \partial Z^{r-1}_{2r-1,-r+2}\right)\right. 
 \ar[r]^(.52){d^r} & Z_{0,r-1}^r\left/\left(Z^{r-1}_{-1,r} + \partial Z^{r-1}_{r-1,1}\right)\right.
}
$$
can be used to show that $d^r$ is surjective. (Note that $Z^r_{r,0}$ is generated by 
$x_{1,r-1} - x_{2,r-2} + \cdots + (-1)^{r-1} x_{r,0}$.)


\begin{example}[Multicomplexes with $\textmc{d}_r \neq 0$ and $d^i = 0$ for all $i \geq 2$]
\end{example}

The preceeding examples show that the spectral sequence associated to a multicomplex 
with $\textmc{d}_r = 0$ for all $r \geq 2$ may not degenerate at $E^2$ (or even $E^r$ where $r$ is
arbitrarily large). These examples can be modified to show that there exist multicomplexes
where $\textmc{d}_r \neq 0$ for $r$ arbitrarily large but the associated spectral sequences degenerate
at $E^2$.

\medskip
We begin with a multicomplex where $\textmc{d}_2 \neq 0$ but its associated spectral sequence
degenerates at $E^2$.  Consider the following first quadrant multicomplex
$$
\xymatrix{
0 \ar[d]^(.4){\textmc{d}_0} & 0 \ar[d]^(.4){\textmc{d}_0} \ar[l]_(.45){\textmc{d}_1} & 0 \ar[d]^(.4){\textmc{d}_0} \ar[l]_(.45){\textmc{d}_1} \\
<x_{0,1}> \ar[d]^(.4){\textmc{d}_0} & <x_{1,1}> \ar[d]^(.4){\textmc{d}_0} \ar[l]_(.45){\textmc{d}_1} & 0 \ar[d]^(.4){\textmc{d}_0} \ar[l]_(.45){\textmc{d}_1} \\
0  & <x_{1,0}> \ar[l]_(.45){\textmc{d}_1} & <x_{2,0}> \ar[l]_(.45){\textmc{d}_1} \ar[llu]|(.37){\textmc{d}_2} |!{[ul];[l]}\hole
}
$$
where the groups $X_{p,q} = 0$ for $p+q > 2$, the homomorphisms $\textmc{d}_0$ and $\textmc{d}_1$
are the same as in Example \ref{doublecomplex}, and $\textmc{d}_2(x_{2,0}) = x_{0,1}$.
The homomorphisms $\textmc{d}_i:X_{p,q} \rightarrow X_{p-i,q+i-1}$ satisfy
$\sum_{i+j = n} \textmc{d}_i \textmc{d}_j = 0 \text{ for all }n$ trivially, and the
assembled chain complex associated to this multicomplex is as follows.
$$
\xymatrix{
\cdots & 0 \ar[r]^(.4){\textmc{d}_0}  \ar[dr]^{\textmc{d}_1} & 0 \\
\cdots & 0 \ar@{}[u]|{\oplus} \ar[r]^(.4){\textmc{d}_0} \ar[dr]^{\textmc{d}_1} & <x_{2,0}> \ar@{}[u]|{\oplus}
\ar[r]^{\textmc{d}_0} \ar[ddr]|(.33){\textmc{d}_2} |!{[d];[dr]}\hole \ar[dr]^{\textmc{d}_1} & 0 & & \\
\cdots & 0 \ar@{}[u]|{\oplus} \ar[r]^(.4){\textmc{d}_0} \ar[dr]^{\textmc{d}_1} & <x_{1,1}> 
\ar@{}[u]|{\oplus} \ar[r]^{\textmc{d}_0} \ar[dr]^{\textmc{d}_1} & <x_{1,0}> \ar@{}[u]|{\oplus} 
\ar[r]^(.63){\textmc{d}_0} \ar[dr]^{\textmc{d}_1} & 0 & \\
\cdots & 0 \ar@{}[u]|{\oplus} \ar[r]^(.4){\textmc{d}_0} & 0 \ar@{}[u]|{\oplus} 
\ar[r]^{\textmc{d}_0} & <x_{0,1}> \ar@{}[u]|{\oplus} \ar[r]^(.63){\textmc{d}_0} & 0 \ar@{}[u]|{\oplus}
\ar[r]^-{\textmc{d}_0} & 0\\
\cdots & 0 \ar@{}[u]|{\|} \ar[r]^(.4){\partial_3} & {(CX)}_2 \ar@{}[u]|{\|} 
\ar[r]^{\partial_2} & {(CX)}_1 \ar@{}[u]|{\|} \ar[r]^(.63){\partial_1} & 0 \ar@{}[u]|{\|}
\ar[r]^-{\partial_0} & 0 \ar@{}[u]|{\|}
}
$$
Referring back to Example \ref{doublecomplex} we see that the diagram
$$
\xymatrix{
Z^r_{2,0} \ar[r]^(.5){\partial} \ar[d] & Z^r_{2-r,r-1} \ar[d]\\
Z^r_{2,0}\left/\left(Z^{r-1}_{1,1} + \partial Z^{r-1}_{r+1,-r+2}\right)\right. 
\ar[r]^(.5){d^r} & Z^r_{2-r,r-1}\left/\left(Z^{r-1}_{1-r,r} + \partial Z^{r-1}_{1,1}\right)\right.
}
$$
shows that $d^r = 0$ for all $r \geq 2$ since $\partial_2(x_{1,1} - x_{2,0}) = 0$,
and hence the associated spectral sequence degenerates at $E^2$.

\medskip
To generalize this example to produce a multicomplex where $\textmc{d}_r \neq 0$ for $r$ arbitrarily
large but the associated spectral sequence still degenerates at $E^2$ we start with the double
complex from Example \ref{double_r} and add a single homomorphism $\textmc{d}_r$ defined by 
$\textmc{d}_r(x_{r,0}) = (-1)^r x_{0,r-1}$.  Further details are left to the reader.


\begin{example}[A multicomplex with $d^2 \neq 0$, $\textmc{d}_2 \neq 0$ and $d^2 \neq \textmc{d}_2$]\label{nonzero}
\end{example}

Consider the following first quadrant multicomplex
$$
\xymatrix{
0 \ar[d]^(.4){\textmc{d}_0} & 0 \ar[d]^(.4){\textmc{d}_0} \ar[l]_(.45){\textmc{d}_1} & 0 \ar[d]^(.4){\textmc{d}_0} \ar[l]_(.55){\textmc{d}_1} \\
<x_{0,1}, \tilde{x}_{0,1}> \ar[d]^(.4){\textmc{d}_0} & <x_{1,1}> \ar[d]^(.4){\textmc{d}_0} \ar[l]_(.45){\textmc{d}_1} & 0 \ar[d]^(.4){\textmc{d}_0} \ar[l]_(.55){\textmc{d}_1} \\
0  & <x_{1,0}> \ar[l]_(.45){\textmc{d}_1} & <x_{2,0}, \tilde{x}_{2,0}> \ar[l]_(.55){\textmc{d}_1} \ar[llu]|(.37){\textmc{d}_2} |!{[ul];[l]}\hole
}
$$
where the groups $X_{p,q} = 0$ for $p+q > 2$, and the homomorphisms 
$\textmc{d}_i$ for $i=0,1,2$ satisfy the following.
$$
\begin{array}{ccr}
\textmc{d}_0(x_{1,1}) & = & x_{1,0}\\
\textmc{d}_1(x_{1,1}) & = & x_{0,1}\\
\textmc{d}_1(x_{2,0}) & = & x_{1,0}\\
\textmc{d}_1(\tilde{x}_{2,0}) & = & 0\\
\textmc{d}_2(\tilde{x}_{2,0}) & = & \tilde{x}_{0,1}\\
\textmc{d}_2(x_{2,0}) & = & 0
\end{array}
$$
The homomorphisms $\textmc{d}_i:X_{p,q} \rightarrow X_{p-i,q+i-1}$ satisfy
$\sum_{i+j = n} \textmc{d}_i \textmc{d}_j = 0 \text{ for all }n$ trivially, and the
assembled chain complex associated to this multicomplex is as follows.
$$
\xymatrix{
\cdots & 0 \ar[r]^(.4){\textmc{d}_0}  \ar[dr]^{\textmc{d}_1} & 0 \\
\cdots & 0 \ar@{}[u]|{\oplus} \ar[r]^(.4){\textmc{d}_0} \ar[dr]^{\textmc{d}_1} & <x_{2,0}, \tilde{x}_{2,0}> \ar@{}[u]|{\oplus} \ar[r]^{\textmc{d}_0} 
\ar[dr]^{\textmc{d}_1} \ar[ddr]|(.33){\textmc{d}_2} |!{[d];[dr]}\hole & 0 & & \\
\cdots & 0 \ar@{}[u]|{\oplus} \ar[r]^(.4){\textmc{d}_0} \ar[dr]^{\textmc{d}_1} & <x_{1,1}> 
\ar@{}[u]|{\oplus} \ar[r]^{\textmc{d}_0} \ar[dr]^{\textmc{d}_1} & <x_{1,0}> \ar@{}[u]|{\oplus} 
\ar[r]^(.62){\textmc{d}_0} \ar[dr]^{\textmc{d}_1} & 0 & \\
\cdots & 0 \ar@{}[u]|{\oplus} \ar[r]^(.4){\textmc{d}_0} & 0 \ar@{}[u]|{\oplus} 
\ar[r]^{\textmc{d}_0} & <x_{0,1}, \tilde{x}_{0,1}> \ar@{}[u]|{\oplus} \ar[r]^(.62){\textmc{d}_0} & 0 \ar@{}[u]|{\oplus}
\ar[r]^-{\textmc{d}_0} & 0\\
\cdots & 0 \ar@{}[u]|{\|} \ar[r]^(.4){\partial_3} & {(CX)}_2 \ar@{}[u]|{\|} 
\ar[r]^{\partial_2} & {(CX)}_1 \ar@{}[u]|{\|} \ar[r]^(.62){\partial_1} & 0 \ar@{}[u]|{\|}
\ar[r]^-{\partial_0} & 0 \ar@{}[u]|{\|}
}
$$
The homology $H_n((CX)_\ast,\partial)$ of the assembled chain complex
is trivial for all $n \in \mathbb{Z}_+$, the $E^1$ term of the associated spectral
sequence is
$$
\xymatrix{
0         & 0 \ar[l]_(.3){\textmc{d}_1} & 0 \ar[l]_(.58){\textmc{d}_1}                    \\
<x_{0,1}, \tilde{x}_{0,1}> & 0 \ar[l]_(.3){\textmc{d}_1} & 0 \ar[l]_(.58){\textmc{d}_1}   \\
0         & 0 \ar[l]_(.3){\textmc{d}_1} & <x_{2,0}, \tilde{x}_{2,0}> \ar[l]_(.58){\textmc{d}_1}
}
$$
where $E^1_{s,t} = 0$ for all $s+t > 2$, and the $E^2$ term is isomorphic to
the $E^1$ term. The image of the homomorphism induced by $\textmc{d}_2$ does not
include the class determined by $x_{0,1}$. However, the differential $d^2$
in the spectral sequence is onto. Therefore, $d^2$ is not the same as the
homomorphism induced by $\textmc{d}_2$.

\bigskip
Note that additional examples can be constructed where the homology is nontrivial
by adding more generators. Examples \ref{doublecomplex}, \ref{double_r}, and \ref{nonzero}
were constructed to have trivial homology in order to make it easy to see that $d^r$ is
surjective. Also, it should be clear at this point how to construct examples where $d^r$ is
not induced from $\textmc{d}_r$ for several different values of $r$: simply
combine the above examples using more (disjoint) generators.


\bigskip\noindent
Acknowledgments:  I would like to thank Jim Stasheff and other members of the University
of Pennsylvania's Deformation Theory Seminar, in particular Tom Hunter and Ron Umble,
for introducing me to multicomplexes and pointing out that the Morse-Bott-Smale chain complex
is in fact a multicomplex.

\nocite{WeiAni}


\def\cprime{$'$}

\end{document}